\documentclass[11pt]{article}

\usepackage{amssymb,amsmath,amsfonts}
\usepackage{color}
\definecolor{darkgreen}{rgb}{0,.3,0}
\definecolor{darkred}{rgb}{.3,0,0}
\definecolor{darkcyan}{rgb}{0,.3,.3}
\definecolor{darkmagenta}{rgb}{.3,0,.3}
\usepackage[dvips,
            colorlinks,%
            linkcolor=darkred,%
            pagecolor=red,%
            filecolor=darkcyan,%
            citecolor=darkgreen,%
            urlcolor=darkmagenta]{hyperref} 
\ifx\pdfoutput\undefined
\usepackage[dvips]{graphicx}
\else
\usepackage[pdftex]{graphicx}
\fi

\def\TheLAYER#1#2#3{{\cal L}_{#2;#3}#1}

\def\demo#1{\par\indent\emph{#1.}}
\def\enddemo{\smallskip\par}
\def\proclaim#1{\par\textbf{#1.}\sl}
\def\endproclaim{\rm\smallskip}
\def\qed{$\blacktriangle$}
\def\wordDefinition{D\,e\,f\,i\,n\,i\,t\,i\,o\,n\,}
\def\wordProposition{P\,r\,o\,p\,o\,s\,i\,t\,i\,o\,n\,\refstepcounter{Claim}}
\def\wordLemma{L\,e\,m\,m\,a\,\refstepcounter{Claim}}
\def\wordTheorem{T\,h\,e\,o\,r\,e\,m\,\refstepcounter{Claim}}
\def\wordCorollary{C\,o\,r\,o\,l\,l\,a\,r\,y\,\refstepcounter{Claim}}
\def\wordRemark{R\,e\,m\,a\,r\,k\,}

\title{Asymptotics for the number of $n$-qua\-si\-groups of order $4$}
\author{V.\,N.\,Potapov and D.\,S.\,Krotov%

\footnote{%
Vladimir Potapov, Denis Krotov,
Sobolev Institute of Mathematics,
4 Acad. Koptyug avenue 630090 Novosibirsk \newline
\texttt{vpotapov@math.nsc.ru}, \texttt{krotov@math.nsc.ru}
\newline\indent
The work of the first author is supported by RFBR (project 05--01--00364).
}
}
\begin{document}
\maketitle
\begin{abstract}
The asymptotic form of the number of $n$-qua\-si\-groups of order $4$ is
$3^{n+1}2^{2^n +1}(1+o(1))$.

\emph{Keywords}: $n$-qua\-si\-groups, MDS codes, decomposability, reducibility.

\emph{MSC}: 20N15, 05B15, 94B25.
\end{abstract}

An algebraic system that consists
of a set
$\Sigma$
of cardinality
$|\Sigma|=k$
and an $n$-ary operation
$f: \Sigma^n\rightarrow \Sigma$
uniquely invertible by
each of its arguments
is called an {\it $n$-qua\-si\-group
of order $k$.}
The function $f$ can also be referred to as an {\it $n$-qua\-si\-group of order $k$} (see \cite{Belousov}).
The value table of an $n$-qua\-si\-group of order $k$
is called a {\it Latin $n$-cube of dimension $k$}
(if $n=2$, a {\it  Latin square}).
Furthermore, there is a one-to-one correspondence
between the $n$-qua\-si\-groups
and the distance $2$ MDS codes of length $n+1$.

It is not difficult to show that for each $n$ there exist only two
$n$-qua\-si\-groups of order $2$ and
$3 \cdot 2^n$
different $n$-qua\-si\-groups of order $3$,
which constitute one equivalence class.
In this work we study properties
of $n$-qua\-si\-groups of order $4$
and derive the asymptotic representation
$3^{n+1}2^{2^n +1}(1+o(1))$
for their number.
The results of the research were announced in \cite{KroPot:Lyap}.
For $k>4$, the asymptotic form of the number of $n$-qua\-si\-groups
and even the asymptotic form of its logarithm remain unknown.

 In Sections~1--4 we give necessary definitions and statements on
 quaternary distance $2$ MDS codes and double-codes (\hyperref[s:1]{Section~1}),
 linear double-codes (\hyperref[s:2]{Section~2}),
 $n$-qua\-si\-groups of order $4$ (\hyperref[s:3]{Section~3}),
 semilinear $n$-qua\-si\-groups of order $4$ (\hyperref[s:4]{Section~4}).
 In \hyperref[s:5]{Section~5} we prove that almost all (as
 $n\to\infty$)
 $n$-qua\-si\-groups of order $4$ are semilinear and establish
 asymptotically tight bounds on their number.

 In addition to the main result, the following lemmas
 can be viewed as stand-alone results:
  \hyperref[l:1]{Lemma~1} on a linear anti-layer in
 a double-MDS-code, \hyperref[l:4]{Lemma~4} on a semilinear layer in
 an $n$-qua\-si\-group, as well as \hyperref[l:2]{Lemmas~2} \hyperref[l:3]{and~3}
 on the decomposability of double-MDS-codes and $n$-qua\-si\-groups,
 proved in \cite{Kro:2002ACCT:double-codes,Kro:2codes}, and their \hyperref[c:3]{Corollary~3}.

\section{MDS codes and double-codes}\label{s:1}

Let
$\Sigma=\{0,1,2,3\}$
and $n$ be a natural number.
In this paper we study subsets of
$\Sigma^n$
and functions defined on
$\Sigma^n$
that
have some properties specified below.
The elements of
$\Sigma^n$
will be called {\it vertices.}
Denote by $[n]$ the set of natural numbers from $1$ to $n$.
Given
$\bar y = (y_1,\dots,y_n)$,
we put
$\bar y^{(i)}\#x = (y_1,\dots,y_{i-1},x,y_{i+1},\dots,y_n)$.

Assume
$\bar x\in \Sigma^n$ and $k\in [n]$.
The set
${\cal E}_k(\bar x) \triangleq \{\bar x^{(k)}\#a:a\in \Sigma\}$
is called a
{\it $k$-edge}.
Two different vertices in
$\Sigma^n$
are called
{\it neighbor}
iff they both belong to some $k$-edge,
i.\,e., differ in only one coordinate.

\demo{\wordDefinition}
A set
$C\subset \Sigma^n$
is called
a {\it distance $2$ {\rm MDS} code} ({\it of length} $n$)
(henceforth simply an {\it {\rm MDS} code}) iff
$|{\cal E}_k(\bar x)\cap C| =1$
for every
$\bar x\in \Sigma^n$
and
$k \in [n]$.
Note that
$|C|=|\Sigma^n|/4 =2^{2n-2}$.
\enddemo

\demo{\wordDefinition}
A set
$S\subset \Sigma^n$
is called
a {\it double-code}
iff
$|{\cal E}_k(\bar x)\cap S| =2$
for every
$\bar x\in S$
and
$k \in [n]$.
\enddemo

\demo{\wordDefinition}
A double-code
$S\subset \Sigma^n$
is called
a {\it double-{\rm MDS}-code}
iff
$|S|=|\Sigma^n|/2 = 2^{2n-1}$.
In other words,  a set
$S\subset \Sigma^n$
is a double-MDS-code iff
$|{\cal E}_k(\bar x)\cap S|=2$
for every
$\bar x\in \Sigma^n$
and
$k \in [n]$.
Obviously,
$\Sigma^n\backslash S$
also is a double-MDS-code in this case.
\enddemo

Denote by
${\Gamma}(S)$
the {\it adjacency graph}
of a double-code
$S\subset \Sigma^n$
with the vertex set $S$ and the edge set
$\{(\bar x,\bar y): \text{$\bar x$, $\bar y$ are
neighbor vertices in $\Sigma^n$}\}$.

\demo{\wordDefinition}
A nonempty double-code
$S\subset \Sigma^n$
is called
{\it prime}
iff $S$ is a subset of a double-MDS-code
$S'\subset \Sigma^n$
and the graph
${\Gamma}(S)$
is connected.
By way of illustration, we list all up to equivalence
nonempty double-codes in
$\Sigma^2$ (Fig.~1).
\enddemo 
\begin{figure}
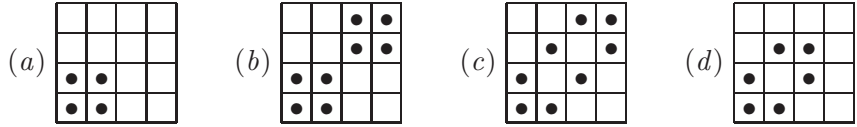

\begin{center}
$
\renewcommand\arraystretch{0.8}
{\mbox{({\em a\/})}\ }
\begin{tabular}[c]{|@{\hspace{0.6ex}}c@{\hspace{0.6ex}}|@{\hspace{0.6ex}}c@{\hspace{0.6ex}}|@{\hspace{0.6ex}}c@{\hspace{0.6ex}}|@{\hspace{0.6ex}}c@{\hspace{0.6ex}}|}
\hline\phantom{\raisebox{-0.2ex}{$\phantom{\bullet}$}}&\phantom{\raisebox{-0.2ex}{$\phantom{\bullet}$}}&\phantom{\raisebox{-0.2ex}{$\phantom{\bullet}$}}&\phantom{\raisebox{-0.2ex}{$\phantom{\bullet}$}}\\
\hline\phantom{\raisebox{-0.2ex}{$\phantom{\bullet}$}}&\phantom{\raisebox{-0.2ex}{$\phantom{\bullet}$}}&\phantom{\raisebox{-0.2ex}{$\phantom{\bullet}$}}&\phantom{\raisebox{-0.2ex}{$\phantom{\bullet}$}}\\
\hline\raisebox{-0.2ex}{$\bullet$}&\raisebox{-0.2ex}{$\bullet$}&\raisebox{-0.2ex}{$\phantom{\bullet}$}&\raisebox{-0.2ex}{$\phantom{\bullet}$}\\
\hline\raisebox{-0.2ex}{$\bullet$}&\raisebox{-0.2ex}{$\bullet$}&\raisebox{-0.2ex}{$\phantom{\bullet}$}&\raisebox{-0.2ex}{$\phantom{\bullet}$}\\
\hline
\end{tabular}
\qquad 
{\mbox{({\em b\/})}\ }
\begin{tabular}[c]{|@{\hspace{0.6ex}}c@{\hspace{0.6ex}}|@{\hspace{0.6ex}}c@{\hspace{0.6ex}}|@{\hspace{0.6ex}}c@{\hspace{0.6ex}}|@{\hspace{0.6ex}}c@{\hspace{0.6ex}}|}
\hline \raisebox{-0.2ex}{$\phantom{\bullet}$}&\raisebox{-0.2ex}{$\phantom{\bullet}$}&\raisebox{-0.2ex}{$\bullet$}&\raisebox{-0.2ex}{$\bullet$}\\
\hline \raisebox{-0.2ex}{$\phantom{\bullet}$}&\raisebox{-0.2ex}{$\phantom{\bullet}$}&\raisebox{-0.2ex}{$\bullet$}&\raisebox{-0.2ex}{$\bullet$}\\
\hline \raisebox{-0.2ex}{$\bullet$}&\raisebox{-0.2ex}{$\bullet$}&\raisebox{-0.2ex}{$\phantom{\bullet}$}&\raisebox{-0.2ex}{$\phantom{\bullet}$}\\
\hline \raisebox{-0.2ex}{$\bullet$}&\raisebox{-0.2ex}{$\bullet$}&\raisebox{-0.2ex}{$\phantom{\bullet}$}&\raisebox{-0.2ex}{$\phantom{\bullet}$}\\
\hline
\end{tabular}
\qquad 
{\mbox{({\em c\/})}\ }
\begin{tabular}[c]{|@{\hspace{0.6ex}}c@{\hspace{0.6ex}}|@{\hspace{0.6ex}}c@{\hspace{0.6ex}}|@{\hspace{0.6ex}}c@{\hspace{0.6ex}}|@{\hspace{0.6ex}}c@{\hspace{0.6ex}}|}
\hline \raisebox{-0.2ex}{$\phantom{\bullet}$}&\raisebox{-0.2ex}{$\phantom{\bullet}$}&\raisebox{-0.2ex}{$\bullet$}&\raisebox{-0.2ex}{$\bullet$}\\
\hline \raisebox{-0.2ex}{$\phantom{\bullet}$}&\raisebox{-0.2ex}{$\bullet$}&\raisebox{-0.2ex}{$\phantom{\bullet}$}&\raisebox{-0.2ex}{$\bullet$}\\
\hline \raisebox{-0.2ex}{$\bullet$}&\raisebox{-0.2ex}{$\phantom{\bullet}$}&\raisebox{-0.2ex}{$\bullet$}&\raisebox{-0.2ex}{$\phantom{\bullet}$}\\
\hline \raisebox{-0.2ex}{$\bullet$}&\raisebox{-0.2ex}{$\bullet$}&\raisebox{-0.2ex}{$\phantom{\bullet}$}&\raisebox{-0.2ex}{$\phantom{\bullet}$}\\
\hline
\end{tabular}
\qquad 
{\mbox{({\em d\/})}\ }
\begin{tabular}[c]{|@{\hspace{0.6ex}}c@{\hspace{0.6ex}}|@{\hspace{0.6ex}}c@{\hspace{0.6ex}}|@{\hspace{0.6ex}}c@{\hspace{0.6ex}}|@{\hspace{0.6ex}}c@{\hspace{0.6ex}}|}
\hline & & & \\
\hline \raisebox{-0.2ex}{$\phantom{\bullet}  $}&\raisebox{-0.2ex}{$\bullet$}&\raisebox{-0.2ex}{$\bullet$}&\raisebox{-0.2ex}{$\phantom{\bullet}$} \\
\hline \raisebox{-0.2ex}{$\bullet$}&\raisebox{-0.2ex}{$\phantom{\bullet}  $}&\raisebox{-0.2ex}{$\bullet$}&\raisebox{-0.2ex}{$\phantom{\bullet}$} \\
\hline \raisebox{-0.2ex}{$\bullet$}&\raisebox{-0.2ex}{$\bullet$}&\raisebox{-0.2ex}{$\phantom{\bullet}  $}&\raisebox{-0.2ex}{$\phantom{\bullet}$} \\
\hline
\end{tabular}
$
 \caption{\label{f:0} The double-codes $(a)$ and $(c)$ are prime;
 $(b)$ and $(c)$ are double-MDS-codes.}
 \end{center}
\end{figure}
\demo{\wordDefinition}
A double-MDS-code $S$
is {\it splittable}
iff
$S=C_1\cup C_2$
where
$C_1$
and
$C_2$
are disjoint MDS codes.
Nonsplittable double-MDS-codes exist in
$\Sigma^n$
starting from
$n=3$.
A double-MDS-code $S$
is splittable if and only if
${\Gamma}(S)$
is a bipartite graph.
\enddemo

\demo{\wordDefinition}
An {\it isotopy},
or
{\it $n$-isotopy},
we call an ordered collection of $n$ permutations
$\theta_i:\Sigma\to \Sigma$,
$i\in [n]$.
Let
$\bar \theta=(\theta_1,\dots ,\theta_n)$
be an isotopy and
$S\subseteq \Sigma^n$.
Put
$\bar \theta S\triangleq\{(\theta_1 x_1,\dots ,\theta_n x_n):(x_1,\dots ,x_n)\in S\}$.
\enddemo

\demo{\wordDefinition}
Sets
$S_1\subseteq \Sigma^n$
and
$S_2\subseteq \Sigma^n$
are called
{\it equivalent}
iff there exist a coordinate permutation
$\tau: [n] \to [n]$
and an $n$-isotopy
$\bar\theta$
such that
$$
\chi_{S_1}(x_1,\dots,x_{n}) \equiv \chi_{\bar\theta S_2} (x_{\tau(1)},\dots,
x_{\tau(n)});
$$
here and in what follows
$\chi_B$
denote the characteristic function of a set $B$.
\enddemo

Obviously, if two double-codes are equivalent then
they have equivalent adjacency graphs;
they are both double-MDS-codes or neither is a double-MDS-code;
they are both splittable or neither is splittable;
and both are prime or neither is prime.

\proclaim{\wordProposition~1}\label{p:1}
Let $S$ be a splittable double-{\rm MDS}-code and
$\gamma$
be the number of the prime double-codes that $S$ includes.
Then the double-code $S$
includes
exactly
$2^{\gamma}$
different {\rm MDS} codes.
\endproclaim

\demo{Proof}
The number of the MDS codes that $S$ includes
equals the number of the ways
of choosing a part of the bipartite graph
${\Gamma}(S)$.
Since in each of the
$\gamma$
connected components of 
${\Gamma}(S)$
the part can be chosen independently, the number of the ways is
$2^{\gamma}$.
\qed\enddemo

\demo{\wordDefinition}
Let
$S\subseteq \Sigma^n$,
$k\in [n]$, and $y\in \Sigma$.
The set
$$
\TheLAYER{S}{k}{y} \triangleq \{(x_1,\dots,x_{k-1},x_{k}, \dots,x_{n-1}):
(x_1,\dots,x_{k-1},y,x_{k},\dots,x_{n-1})\in S\}
$$
is called the $y$th
{\it layer}
of $S$ in direction $k$.
\enddemo

\proclaim{\wordProposition~2}\label{p:2}
Let
$S,S'\subseteq \Sigma^n$
be some sets,
$k\in [n]$,
and
$\{a,b,c,d\}= \Sigma$.

{\rm (a)}
If $S$ is a double-code $($splittable double-code, double-{\rm MDS}-code$)$, then
$\TheLAYER{S}{k}{a}$
also is a double-code $($splittable double-code, double-{\rm MDS}-code$)$ in
$\Sigma^{n-1}$.

{\rm (b)}
If
$k<k'\in [n]$,
then
$\TheLAYER{(\TheLAYER{S}{k'}{a})}{k}b = \TheLAYER{(\TheLAYER{S}{k}{b})}{k'-1}{a}$.

{\rm (c)}
$\TheLAYER{(S \cap S')}{k}{a} =  \TheLAYER{S}{k}{a} \cap \TheLAYER{S'}{k}{a}$.

{\rm (d)}
If $S$ and $S'$ are double-codes and
$\TheLAYER{S}{k}{a} =\TheLAYER{S'}{k}{a}$,
$\TheLAYER{S}{k}{b} =\TheLAYER{S'}{k}{b}$,
$\TheLAYER{S}{k}{c} =\TheLAYER{S'}{k}{c}$,
then
$\TheLAYER{S}{k}{d} =\TheLAYER{S'}{k}{d}$.

{\rm (e)}
If $S$ is a double-{\rm MDS}-code and
$\TheLAYER{S}{k}{a} =\TheLAYER{S}{k}{b}$,
then
$\TheLAYER{S}{k}{c} =\TheLAYER{S}{k}{d} = \Sigma^{n-1}\backslash \TheLAYER{S}{k}{a}$.
\endproclaim

Let us show that a double-MDS-code is completely defined by any of its nonempty
subsets that are double-codes.

\proclaim{\wordProposition~3 \rm(on unique extension of a double-code)}\label{p:3}
Let
$S_1,S_2\subset \Sigma^n$ be
double-{\rm MDS}-codes. Then

{\rm (a)}
if
$S_0 \subseteq S_1\cap S_2$
is a nonempty double-code, then
$S_1 =S_2;$

{\rm (b)}
if
$S_0 \subseteq S_1\backslash S_2 $
is a nonempty double-code, then
$S_1 =\Sigma^n\backslash S_2$.
\endproclaim

\demo{Proof}
We will prove (a) by induction on $n$.
For
$n=1$
the claim is trivial.
Assume that (a) holds for
$n=m-1$;
let us show that it holds for
$n=m$.
By \hyperref[p:2]{Proposition~2(a)}, we have:
$\TheLAYER{S_0}{1}{a}$
is a double-code,
$\TheLAYER{S_1}{1}{a}$
and
$\TheLAYER{S_2}{1}{a}$
are double-MDS-codes for each
$a\in \Sigma$.
By \hyperref[p:2]{Proposition~2(c)},
$\TheLAYER{S_0}{1}{a}\subseteq \TheLAYER{S_1}{1}{a} \cap \TheLAYER{S_2}{1}{a}$.
Then, by the inductive assumption,
$\TheLAYER{S_1}{1}{a} =\TheLAYER{S_2}{1}{a}$
for all
$a\in \Sigma$
such that
$\TheLAYER{S_0}{1}{a}$
is not empty. By the definition of a double-code, at least two of
the four sets
$\TheLAYER{S_0}{1}{a}$,
$a\in\Sigma$
are nonempty. If there are three nonempty sets, then the equality
$S_1=S_2$
follows from \hyperref[p:2]{Proposition~2(d)}.
Assume that two sets, say
$\TheLAYER{S_0}{1}{2}$
and
$\TheLAYER{S_0}{1}{3}$,
are empty. Then
$\TheLAYER{S_0}{1}{0}=\TheLAYER{S_0}{1}{1}$,
because
$|{\cal E}_1(\bar x)\cap S_0| = 2$
for all
$\bar x\in S_0$.
Hence
$\TheLAYER{S_1}{1}{0}=\TheLAYER{S_1}{1}{1} =\TheLAYER{S_2}{1}{0}=\TheLAYER{S_2}{1}{1}$,
by the inductive assumption.
Then, by \hyperref[p:2]{Proposition~2(e)}, we get
$S_1=S_2$.

(b) Consider
$S'_2 \triangleq \Sigma^n\backslash S_2$.
Since
$S'_2$
is a double-MDS-code and
$S_0\subseteq S_1\cap S_2'$,
it follows from (a) that
$S_2'=S_1$.
\qed
\enddemo

\section{Linear double-codes}\label{s:2}

\demo{\wordDefinition}
A nonempty double-code
$S\subset \Sigma^n$
is called
{\it linear}
iff

\begin{equation}
\chi_{S}(x_1,\dots,x_{n}) \equiv \chi_{S_1}(x_1)\oplus \chi_{S_2}(x_2)\oplus
\dots \oplus \chi_{S_n}(x_n)
\label{e:1}
\end{equation}
where
$S_i$
($1\leq i\leq n$)
are subsets of
$\Sigma$
and
$\oplus$
is the modulo $2$ addition. Obviously,
$S_i$
are double-MDS-codes in
$\Sigma$.
A linear double-code in
$\Sigma^2$
is illustrated in {\rm Fig.~1$(b)$}.
\enddemo

In the following two propositions,
some elementary properties of linear $2$-codes are proved.

\proclaim{\wordProposition~4 \rm(properties of the class of linear double-codes)}\label{p:4}
{\rm (a)}
 The linear double-codes constitute an equivalence class.

{\rm (b)}
 A linear double-code is a splittable double-{\rm MDS}-code.

{\rm (c)}
 The complement of a linear double-code is a linear double-code.

{\rm (d)}
 A double-code $S$ is linear if and only if
 there exist a prime double-code
 $S_0\subset S$
 equivalent to
 $\{0,1\}^n$.

{\rm (e)}
 A linear double-code is uniquely defined
 by the subset of all its vertices of type
 $\bar 0^{(i)}\#y$,
 $i\in [n]$,
 $y\in \Sigma$.

{\rm (f)}
 The number of linear double-codes in
 $\Sigma^n$
 is
 $2\cdot 3^n$.
\endproclaim

\demo{Proof}
The properties (a)--(c) follow from definitions.

\smallskip
(d) {\it Necessity}. By (a), we can assume without loss of generality that
$\chi_{S}(x_1,\dots,x_{n}) \equiv \bigoplus\limits_{i=1}^n \chi_{\{2,3\}}(x_i)$.
In this case
$S_0 \triangleq \{2,3\} \times \{0,1\}^{n-1}$
is a subset of $S$.

\smallskip
{\it Sufficiency}. Suppose that a double-code
$S_0\subset S$
is equivalent to
$\{0,1\}^n$.
Without loss of generality assume
$S_0=\{2,3\} \times \{0,1\}^{n-1}$.
Then
$S_0$
is a subset of the linear double-code
$S'$
where
$\chi_{S'}(x_1,\dots ,x_n) \equiv \bigoplus\limits_{i=1}^n \chi_{\{2,3\}}(x_i)$.
By \hyperref[p:3]{Proposition~3(a)}, we have
$S=S'$.

(e) Indeed, let a double-code $S$ be represented as in \hyperref[e:1]{(1)}.
Put
$\chi^0 \triangleq \chi_S(\bar 0)$
and
$\chi^i(y) \triangleq \chi_S(\bar 0^{(i)}\# y)$, $i\in [n]$;
then we have
\begin{equation}
\chi_S(x_1,\ldots,x_n)\equiv \chi^0 \oplus
\bigoplus_{i=1}^n(\chi^i(x_i)\oplus \chi^0),
\label{e:2}
\end{equation}
which can be easily checked applying
the formula \hyperref[e:1]{(1)} for $\chi_S$.

(f) follows from the representation \hyperref[e:2]{(2)}. Indeed,
we can choose $\chi^0$ in two ways;
then each of the functions
$\chi^i$,
$i\in [n]$,
can be chosen in three ways,
taking into account that $\chi^i$ is the characteristic
function of a double-MDS-code in
$\Sigma$
and
$\chi^i(\bar 0)=\chi^0$.~
\qed\enddemo

The set
$\{0,1\}^n$
(as well as the graph
${\Gamma}( \{0,1\}^n)$)
is called the
{\it Boolean $n$-cube}.
The next proposition follows from definitions and \hyperref[p:2]{Proposition~2}.

\proclaim{\wordProposition~5 \rm(on heritable properties of linear double-codes)}\label{p:5}
{\rm (a)}
 If
 $S\subset \Sigma^n$
 is a linear double-code, then
 $\TheLAYER{S}{k}{y}$
 is a linear double-code.

 {\rm (b)}
 Let
 $S\subset \Sigma^n$
 be a double-code.
 If two layers of $S$ by some direction are linear and coincides,
 then $S$ is a linear double-code.
\endproclaim

The main result of this section is the following lemma,
presenting a partial inversion of p.\,(a) and a partial
strengthening of p.\,(b) of \hyperref[p:5]{Proposition~5}. The lemma claim that
the existence of a linear layer in a splittable double-MDS-code implies
the existence of a layer (``anti-layer'') in the same direction
that complements the former.

\proclaim{\wordLemma~1 \rm (on a linear anti-layer)}\label{l:1}
Let
$S\subset \Sigma^n$
be a splittable double-{\rm MDS}-code and
$L \triangleq \TheLAYER{S}{k}{a}$
be a linear double-code for some
$k\in [n]$
and
$a\in \Sigma$.
Then

{\rm (a)}
 there is
 $b\in \Sigma$
 such that
 $\TheLAYER{S}{k}{b} =\Sigma^{n-1}\backslash L;$

{\rm (b)}
 $\Sigma^{n}\backslash S$
 is a splittable double-{\rm MDS}-code.
\endproclaim

Before proving \hyperref[l:1]{Lemma~1} we introduce the notation
$\neg(\alpha_1,\alpha_2,\dots,\alpha_n) \triangleq
(\alpha_1\oplus1,\alpha_2\oplus1,\dots,\alpha_n\oplus1)$
where
$\alpha_i\in \{0,1\}$,
and prove two auxiliary propositions.

\proclaim{\wordProposition~6}\label{p:6}
Let
$\{P_1,P_2,P_3\}$
be a partition of the Boolean $n$-cube with
$n\geq 4$
into three nonempty sets{\rm :}
$P_1 \cup P_2 \cup P_3 = \{0,1\}^n$.
And assume the following holds{\rm :}

$(*)$ for every
$k\in [n]$
and every
$b \in \{0,1\}$
at least one set {\rm (}layer{\rm )} of
$\TheLAYER{P_1}{k}{b}$,
$\TheLAYER{P_2}{k}{b}$,
$\TheLAYER{P_3}{k}{b}$
is empty.

Then
$\{P_1,P_2,P_3\} =
 \{\{\bar\alpha\},\{\neg\bar\alpha\},\{0,1\}^n\backslash
 \{\bar\alpha,\neg\bar\alpha\} \}$
where
$\bar\alpha\in\{0,1\}^n$.
\endproclaim

\demo{Proof}
Denote by
$N_i\subseteq [n]$
the set of coordinates $k$ whose values are not fixed in
$P_i$,
i.\,e.,
$\TheLAYER{P_i}{k}{0}\neq \emptyset$
and
$\TheLAYER{P_i}{k}{1}\neq \emptyset$.
It is easy to see that the sets
$N_1, N_2, N_3$
are pairwise disjoint (if, for example,
$k\in N_1\cap N_2$,
then $(*)$ implies
$\TheLAYER{P_3}{k}{0} = \emptyset$
and
$\TheLAYER{P_3}{k}{1} = \emptyset$,
which contradicts the nonemptiness of
$P_3$).
So, the obvious relation
$2^n= |P_1|+|P_2|+|P_3| \leq 2^{|N_1|}+2^{|N_2|}+2^{|N_3|}$
yields
$\{N_1,N_2,N_3\}=\{\emptyset,\emptyset, [n]\}$
and
$\{P_1,P_2,P_3\} =
 \{\{\bar\alpha\},\{\bar\beta\},\{0,1\}^n\backslash
 \{\bar\alpha,\bar\beta\} \}$.
The hypothesis $(*)$ implies that
$\bar\beta=\neg\bar\alpha$.
\qed\enddemo

\proclaim{\wordProposition~7}\label{p:7}
Let $S$
be a double-{\rm MDS}-code in
$\Sigma^n$,
$n\geq 3$,
and
$k\in [n]$.
Let
$P_0$,
$P_1$,
$P_2$,
$P_3$
be the intersections of the four layers of $S$ in direction $k$ with the Boolean
$(n-1)$-cube, i.\,e.,
$P_i \triangleq \TheLAYER{S}{k}{i} \cap \{0,1\}^{n-1}$.
Assume that at least one of the following holds{\rm :}

{\rm (a)}
 $n=3$,
 $P_i=\{0,1\}^2$
 for some $i$, and
 $P_i\neq \emptyset$
 for all
 $i\in\{0,1,2,3\};$

{\rm (b)}
 $\{P_0, P_1, P_2, P_3\} = \{\{0,1\}^{n-1}, \{\bar\alpha\}, \{\bar\beta\},
  \{0,1\}^{n-1} \backslash \{\bar\alpha,\bar\beta\}\}$
 where
 $\bar\alpha\in \{0,1\}^{n-1}$
 and
 $\bar\beta=\neg\bar\alpha$.

Then the double-codes $S$ and
$\Sigma^n \backslash S$
are unsplittable.
\endproclaim

\demo{Proof}
(a) There are two nonequivalent cases for a choice of the sets
$P_i$.
It is not difficult to check 
(we leave this to the reader)
that in each case an attempt to recover the double-MDS-code $S$
leads to an unsplittable double-MDS-code with the unsplittable complement.

(b)
Without loss of generality we can assume that
$k=n$,
$\bar\alpha = 0^{n-1}$, $\bar\beta = 1^{n-1}$,
$$
P_0=\{0,1\}^{n-1},\quad P_1=\{\bar\alpha\},\quad  P_2=\{\bar\beta\},
\quad  P_3=\{0,1\}^{n-1}
\backslash \{\bar\alpha,\bar\beta\}
$$
(otherwise we can select a suitable coordinate permutation and
isotopy and consider an equivalent double-code
that satisfies this assumption).
We will argue by induction on $n$.
The basis of induction, the case of
$n=3$,
is considered in p.\,(a).
Assume that the statement holds for
$n=m-1$.
Let us show that it holds for
$n=m$
as well. Consider the intersections of the layers
$\TheLAYER{S}{k}{0}$,
$\TheLAYER{S}{k}{1}$,
$\TheLAYER{S}{k}{2}$,
$\TheLAYER{S}{k}{3}$
with the set
$E\triangleq\{2,3\}\times \{0,1\}^{n-2}$,
which is equivalent to the Boolean
$(n-1)$-cube
$\{0,1\}^{n-1}$
and is a ``neighbor cube'' to it:
$$
Q_i \triangleq \{2,3\}\times \{0,1\}^{n-2} \cap \TheLAYER{S}{1}{i}.
$$
Fig.~2 illustrates the situation.
\begin{figure}
 \begin{center}
 \includegraphics[scale=0.5]{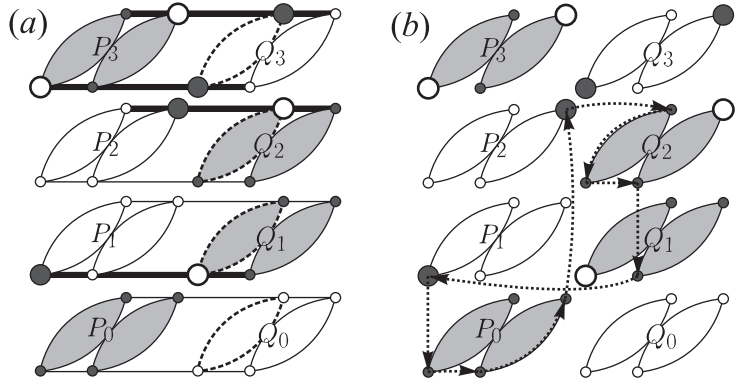}
 \caption{\label{f:1} An illustration for \hyperref[p:7]{Proposition~7}.}
 \end{center}
\end{figure}

$(*)$ We claim that the sets
$Q_0$,
$Q_1$,
$Q_2$,
$Q_3$
are defined up to four elements. More exactly,
\begin{equation}
Q_0=\emptyset,\quad  Q_1=E \backslash \{\bar\alpha'\},\quad  Q_2=E \backslash \{\bar\beta'\},
\quad  Q_3=\{\bar\alpha'',\bar\beta''\}
\label{e:3}
\end{equation}
where
$\bar\alpha',\bar\alpha''\in \{(2,0,\dots ,0),(3,0,\dots ,0)\}$
and
$\bar\beta',\bar\beta''\in \{(2,1,\dots ,1),(3,1,\dots ,1)\}$.
Indeed, the set
$\{0,1\}^{n-1}\cup E$
can be split into
the $1$-edges of type
${\cal E}_1(\bar x)$, $\bar x\in \{0\}\times\{0,1\}^{n-2}$.
Since $S$ is a double-MDS-code,
every such $1$-edge contains two vertices from
$P_i\cup Q_i$
for each
$i\in\{0,1,2,3\}$.
In particular,

$\bullet$
 if such $1$-edge contains two vertices from
 $P_i$,
 then it does not contain vertices from
 $Q_i$;

$\bullet$
 if it does not contain vertices from
 $P_i$,
 then it contains two vertices from
 $Q_i$.

According to \hyperref[e:3]{(3)} these two rules define all vertices of
$Q_i$,
$i=0,1,2,3$,
except for the four cases (Fig.~2$(a)$, the bold horizontal lines):

$\bullet$
 the $1$-edge
 ${\cal E}_1(0,0,\dots ,0)$
 contains exactly one vertex
 $(0,0,\dots ,0)$
 from
 $P_1$,

$\bullet$
 the $1$-edge
 ${\cal E}_1(0,0,\dots ,0)$
 contains exactly one vertex
 $(1,0,\dots ,0)$
 from
 $P_3$,

$\bullet$
 the $1$-edge
 ${\cal E}_1(0,1,\dots ,1)$
 contains exactly one vertex
 $(1,1,\dots ,1)$
 from
 $P_2$,

$\bullet$
 the $1$-edge
 ${\cal E}_1(0,1,\dots ,1)$
 contains exactly one vertex
 $(0,1,\dots ,1)$
 from
 $P_3$.

In each of the cases we have a choice of a vertex of
$Q_i$
for the respective $i$. This choice corresponds to the choice of
$\alpha',\alpha'',\beta',\beta''$.
The claim $(*)$ is proved.

Since $S$ is a double-MDS-code,
every vertex from $E$ belongs to exactly two sets
$Q_i$.
So, it follows directly from \hyperref[e:3]{(3)} that
$\bar\alpha'=\bar\alpha''$ and $\bar\beta'=\bar\beta''$.
Without loss of generality we can assume that
$\bar\alpha'=\bar\alpha''=(2,0,\dots ,0)$.
Thus, it suffices to consider the two cases:
$\bar\beta'=\bar\beta''=(2,1,\dots ,1)$
(Fig.~2$(a)$) and
$\bar\beta'=\bar\beta''=(3,1,\dots ,1)$
(Fig.~2$(b)$).

\smallskip
1. {\em Case}
$\bar\beta'=\bar\beta''=(2,1,\dots ,1)$
(Fig.~2$(a)$).
In this case we can use the inductive assumption.
Indeed, consider the set
$\Sigma^{n-1}\backslash \TheLAYER{S}{1}{2}$.
Its layers in the last direction intersected with the Boolean
$(n-2)$-cube
coincide with
$\{0,1\}^{n-2}$,
$\{(0,\dots ,0)\}$,
$\{(1,\dots ,1)\}$,
and
$\{0,1\}^{n-1} \backslash \{(0,\dots ,0),(1,\dots ,1)\}$
(see Fig.~2$(a)$, the dotted lines).
By the inductive assumption, the double-codes
$\Sigma^{n-1}\backslash \TheLAYER{S}{1}{2}$
and
$\TheLAYER{S}{1}{2}$
are unsplittable. Hence,
$\Sigma^{n}\backslash S$
and
$S$ are unsplittable.

\smallskip
2. {\em Case}
$\bar\beta'=\bar\beta''=(3,1,\dots ,1)$
(Fig.~2$(b)$).
In this case we can find a cyclic path
of odd length $2n+3$ in ${\Gamma}(S)$:
$$
\gathered
(0000... 00,~\underbrace{1000... 00,~1100... 00,~1110... 00,~\cdots,~1111... 10}_{n-1},~1111... 12,\phantom{~3000... 01,~0000... 01)}
\\
\phantom{(0000... 00,~}\underbrace{2111... 12,~2011... 12,~2001... 12,~\cdots,~2000... 02}_{n-1},~3000... 02,~3000... 01,~0000... 01)
\endgathered
$$
(Fig.~2$(b)$, the dotted lines); this implies that the graph
is not bipartite and the double-code $S$ is unsplittable
by definition. Similarly, the odd cyclic path
$$
\gathered
(2000... 00,~{3000... 00,~3100... 00,~3110... 00,~\cdots,~3111... 10},~3111... 12,\phantom{~1000... 01,~2000... 01)}
\\
\phantom{(2000... 00,~}{0111... 12,~0011... 12,~0001... 12,~\cdots,~0000... 02},~1000... 02,~1000... 01,~2000... 01)
\endgathered
$$
in
${\Gamma}(\Sigma^n \backslash S)$
shows that the double-code
$\Sigma^n \backslash S$ is unsplittable.
\qed\enddemo

\demo{Proof of \hyperref[l:1]{Lemma~1}}
(a)
We prove the claim by induction.
The base of induction, the case of
$n=2$,
is trivial.
Assume that the lemma holds for
$n=m-1$.
Let us show that it holds for
$n=m\geq 3$.

By \hyperref[p:4]{Proposition~4(d)} and because the splittability and linearity of
a double-code are preserved under isotopy and coordinate
permutation,
without loss of generality we can assume
$k=n$,
$a=0$,
and the linear double-code $L$ includes
$\{0,1\}^{n-1}$.
Let the sets
$P_0$,
$P_1$,
$P_2$,
and
$P_3$
be defined as in \hyperref[p:7]{Proposition~7}, i.\,e.,
$P_i \triangleq \{0,1\}^{n-1} \cap \TheLAYER{S}{n}{i}$.

It is enough to show that at list one of the sets
$P_1$,
$P_2$,
$P_3$
is empty. Then by \hyperref[p:3]{Proposition~3(b)}
the corresponding layer of $S$ will be the complement of $L$.

$(*)$
Assume the contrary, i.\,e., that each of the sets
$P_1$,
$P_2$,
and
$P_3$
is nonempty.

$(**)$
Then we claim that the sets
$P_1$,
$P_2$,
and
$P_3$
satisfy the hypothesis of \hyperref[p:6]{Proposition~6}.
Since $S$ is a double-MDS-code,
its layers in the given direction constitute a twofold
covering of 
$\Sigma^{n-1}$;
and the sets
$P_0$,
$P_1$,
$P_2$,
and
$P_3$
constitute a twofold covering of 
$\{0,1\}^{n-1}$.
Since
$P_0=\{0,1\}^{n-1}$,
we get that
$P_1$,
$P_2$,
and
$P_3$
are pairwise
disjoint and
$P_1 \cup P_2 \cup P_3=\{0,1\}^{n-1}$.
It is remains to show that for every
$r\in [n-1]$
and
$b \in \{0,1\}$
at least one set of
$\TheLAYER{P_1}{r}{b}$,
$\TheLAYER{P_2}{r}{b}$,
$\TheLAYER{P_3}{r}{b}$
is empty. This fact follows from the inductive assumption.
Indeed, the double-code
$\TheLAYER{S}{r}{b}$
fully satisfies the hypothesis of the lemma, and, by the inductive assumption,
it has a layer
$\TheLAYER{\TheLAYER{S}{r}{b}}{n-1}{i}$,
$i\in\{1,2,3\}$
complementary to the ``linear'' layer
$\TheLAYER{\TheLAYER{S}{r}{b}}{n-1}{0}$.
Using \hyperref[p:2]{Proposition~2(b),(d)}
and the inclusion
$\TheLAYER{\TheLAYER{S}{r}{b}}{n-1}{0} \supset \{0,1\}^{n-2}$,
we get
\begin{eqnarray*}
\TheLAYER{P_i}{r}{b}=\TheLAYER{\bigl(\{0,1\}^{n-1} \cap \TheLAYER{S}{n}{i}\bigr)}{r}{b} =
\{0,1\}^{n-2}\cap \TheLAYER{\TheLAYER{S}{n}{i}}{r}{b}
\\
= \{0,1\}^{n-2}\cap
\TheLAYER{\TheLAYER{S}{r}{b}}{n-1}{i} = \emptyset.
\end{eqnarray*}
The claim $(**)$ is proved.

By \hyperref[p:6]{Proposition~6}, the set $S$
satisfies the hypothesis of \hyperref[p:7]{Proposition~7}.
This means the double-code $S$ is unsplittable,
which contradicts to the hypothesis of the lemma. Thus the assumption $(*)$
is not true, and one of the sets
$P_1$,
$P_2$,
and
$P_3$
is empty.

Suppose
$P_j=\emptyset$.
Then
$b=j$,
$\{0,1\}^{n-1} \subset L \backslash \TheLAYER{S}{n}{b}$;
therefore
$\TheLAYER{S}{n}{b} = \Sigma^{n-1} \backslash L$
by \hyperref[p:3]{Proposition~3(b)}.
The claim (a) of the lemma is proved.

(b) As shown in p.\,(a),
two layers of the double-MDS-code $S$ in direction $k$ are complements to each other (with respect to
$\Sigma^{n-1}$).
The definition of a double-code implies that
the other two layers also are complements to each other. Hence
an appropriate permutation of the layers converts $S$ to its complement
$\Sigma^n \backslash S$
and the splittability of the former means the splittability of the later.
\qed\enddemo

Examples show that the layer linearity hypothesis in \hyperref[l:1]{Lemma~1}
is essential for the existence of a layer complementary to a given one
in a splittable double-MDS-code.

\section{ MDS codes and $n$-qua\-si\-groups}\label{s:3}

\demo{\wordDefinition}
Let
$G\subseteq \Sigma^n = \{0,1,2,3\}^n$;
a function
$f: G \rightarrow \Sigma$
is called
a {\it partial  $n$-qua\-si\-group of order $4$}
iff the equation
\begin{equation}
f(\bar a^{(i)}\#x) =f(a_1,\dots a_{i-1},x,a_{i+1},\dots a_n)= b
\label{e:4}
\end{equation}
has at most one solution
$x\in \Sigma$
for every
$\bar a\in \Sigma^n$
and
$b \in \Sigma$.
If, in addition,
$G=\Sigma^n$,
then the function $f$ is called an {\it $n$-qua\-si\-group} of order $4$
(in what follows we omit the words ``of order $4$'').
In this case the equation \hyperref[e:4]{(4)} has exactly one solution for every
$\bar a\in \Sigma^n$
and
$b \in \Sigma$.
By
$f^{\langle i\rangle}$
we denote
the {\it inversion}
of the $n$-qua\-si\-group $f$ in $i$th argument,
which is defined by the relation
$$
f^{\langle i\rangle}(\bar x)=b \quad \Longleftrightarrow \quad
f(\bar x^{(i)}\#b)=x_i.
$$
Obviously, the inversion of an $n$-qua\-si\-group $f$ in
each argument also is an $n$-qua\-si\-group.
\enddemo

\demo{\wordDefinition}
An $n$-qua\-si\-group
$g :\Sigma^n \rightarrow \Sigma$
is called
an {\it extension}
of a partial
$n$-qua\-si\-group
$f:G\to \Sigma$
iff
$f= g|_G$.
A partial $n$-qua\-si\-group
that have at least one extension is called
{\it extendable}.
\enddemo

\demo{\wordDefinition}
An $n$-qua\-si\-group $f$ is called
{\it reduced}
iff
$f(\bar 0^{(i)}\#a) =a$
for every
$i\in [n]$
and
$a\in \Sigma$.
A permutation
$\tau:\Sigma \to \Sigma$
is called
{\it reduced}
iff
$\tau(0)=0$.
\enddemo

\demo{\wordDefinition}
An $n$-qua\-si\-group $f$ is called
{\it decomposable}
iff there exist an integer $m$,
$2\leq m < n$,
an $(n-m+1)$-qua\-si\-group
$h$,
an $m$-qua\-si\-group  $g$,
and a permutation
$\sigma: [n] \to [n]$
such that
$$
f(x_1,\dots,x_{n}) \equiv h\left(g(x_{\sigma(1)},\dots, x_{\sigma(m)}),
 x_{\sigma(m+1)},\dots, x_{\sigma(n)}\right).
$$
Fig.~3 shows examples of decomposable (a), (c) and indecomposable (b) $3$-qua\-si\-groups.
\enddemo
\begin{figure}
 \begin{center}
 \includegraphics[scale=0.55]{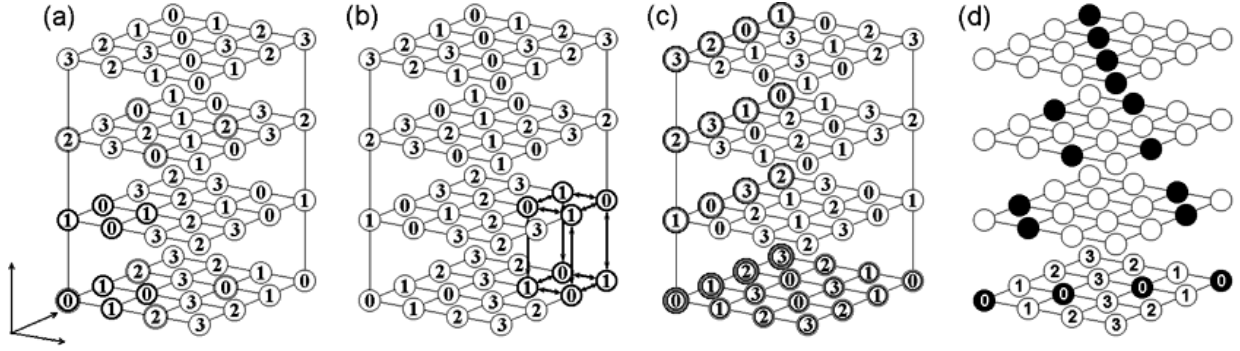}
 \caption{\label{f:3} Examples of $n$-quasigroups:
 (a) a linear $3$-qua\-si\-group, (b) a semilinear indecomposable $3$-qua\-si\-group,
 (c) a nonsemilinear decomposable $3$-quasigroup; (d) a linear $2$-qua\-si\-group and
 the corresponding MDS code.}
 \end{center}
\end{figure}
Take
$f:\Sigma^n\rightarrow \Sigma$
and define the sets
$$
C(f) \triangleq \{(\bar x,f(\bar x)): \bar x\in \Sigma^n\},
\
C_a(f) \triangleq \{\bar x\in \Sigma^n: f(\bar x)=a\},
\
S_{a,b}(f) \triangleq C_a(f)\cup C_b(f).
$$
It follows from definitions that

\proclaim{\wordProposition~8}\label{p:8}
{\rm (a)}
 The mapping
 $C(\cdot)$
 is a one-to-one correspondence between
 the set of all $n$-qua\-si\-groups and the set of all {\rm MDS} codes of length
 $n+1$ $($see Fig.~$3(d)$ for example$)$.

{\rm (b)}
 A function
 $f: \Sigma^n \to \Sigma$
 is an $n$-qua\-si\-group if and only if
 the sets
 $C_a(f)$
 are pairwise disjoint {\rm MDS}-codes for all
 $a\in \Sigma$.

{\rm (c)}
 A function
 $f: \Sigma^n \to \Sigma$
 is an $n$-qua\-si\-group if and only if for every
 different $a$ and $b$ in
 $\Sigma$
 the set
 $S_{a,b}(f)$
 is a splittable double-{\rm MDS}-code.
\endproclaim

\demo{\wordDefinition}
$n$-qua\-si\-groups $f$ and $g$ are called
{\it equivalent}
iff there exist a permutation
$\sigma: [n] \to [n]$
and
an $(n+1)$-isotopy
$\bar \tau = (\tau_0,\tau_1,\dots,\tau_{n})$
such that
$$
f(x_1,\dots,x_{n})
\equiv \tau_0g(\tau_1x_{\sigma(1)},\dots,\tau_{n}x_{\sigma(n)}).
$$
A set of $n$-qua\-si\-groups is called
{\it closed under equivalence}
iff it contains $n$-qua\-si\-groups
together with their equivalence classes.
\enddemo

It follows from definitions that if $n$-qua\-si\-groups $f$ and $g$
are equivalent, then the MDS codes
$C(f)$
and
$C(g)$
are equivalent too. Moreover, an $n$-qua\-si\-group $f$ and its inversion
$f^{\langle i\rangle}$,
$i\in [n]$
correspond to the equivalent MDS codes
$C(f)$
and
$C(f^{\langle i\rangle})$.
For
$n\geq 3$,
there are examples in which an $n$-qua\-si\-group and its inversion are not equivalent.
Thus the equivalence of MDS codes does not imply that the corresponding
$n$-qua\-si\-groups are equivalent.
However, we easily see

\proclaim{\wordProposition~9}\label{p:9}
{\rm(a)}
 Equivalent $n$-qua\-si\-groups are decomposable or
nondecomposable simultaneously.

{\rm(b)}
 If $n$-qua\-si\-group $f$ decomposable, then so are its inversions
 $f^{\langle i\rangle}$,
 $i\in [n]$.
\endproclaim

\proclaim{\wordProposition~10}\label{p:10}
Let
$f:\Sigma^n\rightarrow \Sigma$
be an $n$-qua\-si\-group. Then there exist a unique isotopy
$(\tau_0,\tau_1,\dots,\tau_n)$
with
$\tau_0 = (0,a)$,
$a\in \Sigma$
and reduced permutations
$\tau_1,\dots,\tau_n:\Sigma \to \Sigma$
 such that
\begin{equation}
f(\bar x) \equiv \tau_0g(\tau_1 x_1,\tau_2 x_2,\dots,\tau_n x_n)
\label{e:5}
\end{equation}
where $g$ is a reduced $n$-qua\-si\-group,
$\bar x = (x_1, x_2, \ldots, x_n)$.
\endproclaim

\demo{Proof}
From \hyperref[e:5]{(5)} we get

\begin{eqnarray}\nonumber
  &&\tau_0(0)=f(0,\dots,0),\qquad\text{i.\,e.,}\quad\tau_0=(0,f(0,\dots,0)),
  \\ \label{e:6}
  &&\tau_i (b)=\tau_0^{-1} f(\bar 0^{(i)}\#{b}),\qquad i=1,\ldots,n,
  \\ \nonumber
  &&g(\bar x) = \tau_0^{-1} f\bigl(\tau_1^{-1}x_1,\tau_2^{-1}x_2, \dots,\tau_n^{-1}x_n\bigr);
\end{eqnarray}
%
this yields the uniqueness of the representation.
On the other hand, it is directly verified that
if we define
$\tau_0,\tau_1,\dots,\tau_n$
and $g$ by the equations \hyperref[e:6]{(6)},
then the conditions of the proposition will be satisfied.
\qed\enddemo

Let
$V_n$
be the set of all $n$-qua\-si\-groups of order $4$. Denote by
$R_n\subseteq V_n$
the set of all decomposable $n$-qua\-si\-groups and by
$V^{\star}_n\subset V_n$
the set of all reduced $n$-qua\-si\-groups.
For an arbitrary subset of
$V_n$
denoted by a capital letter with index, for example
$W_n$,
we introduce the following notation:
$W^{\star}_n\triangleq W_n\cap V^{\star}_n$,
$w_n\triangleq|W_n|$,
and
$w^{\star}_n\triangleq|W^{\star}_n|$.

It follows directly from \hyperref[p:10]{Proposition~10} that

\proclaim{\wordCorollary~1}\label{c:1}
Let
$W_n\subseteq V_n$
be a set of $n$-qua\-si\-groups of order $4$ closed under equivalence. Then
$w_n =4\cdot 6^n w^{\star}_n$.
\endproclaim

A partial $n$-qua\-si\-group
$g:G\to \Sigma$ is called
{\it compatible}
with an $n$-qua\-si\-group $f$ iff
$f(\bar x)\neq g(\bar x)$
for every
$\bar x$
from $G$. Denote by
$F(g)$
the set of all $n$-qua\-si\-groups compatible with an $n$-qua\-si\-group $g$.

\proclaim{\wordProposition~11}\label{p:11}
Let $g$ be
an $n$-qua\-si\-group,
$W_n\subseteq V_n$
be a set of $n$-qua\-si\-groups that is closed under equivalence. Then
$|F(g)\cap W_n| \leq 3^{n+1}w^*_{n}$.
\endproclaim

\demo{Proof}
Consider the set
$T\subset\Sigma^{n}$
that consists of the vertices differing from
$(0,\dots,0)\in \Sigma^{n}$
in at most one position.
Given partial $n$-qua\-si\-group
$t:T\rightarrow\Sigma$,
consider the set
$W_{n}(t)$
of its extensions from the class
$W_{n}$, i.\,e.,
$W_{n}(t) \triangleq \{f\in W_{n}: f|_T=t\}$.
Since
$W_{n}$
is closed under equivalence, we have
$|W_{n}(t)| = w^{\star}_{n}$.

It is easy to see that there are exactly
$3^{n+1}$
different partial $n$-qua\-si\-groups
$t:T\rightarrow \Sigma$
compatible with a given $n$-qua\-si\-group $g$.
Since an $n$-qua\-si\-group
$f\in W_n(t)$
is compatible with $g$ only if
$t=f|_T$
is compatible with $g$, the number of the $n$-qua\-si\-groups from
$W_n$ that are compatible with $g$ does not exceed
$3^{n+1}w^\star_n$.
\qed\enddemo

Let
$q:\Sigma^{n-1}\times A\to \Sigma$
be a partial $n$-qua\-si\-group,
$A \subseteq\Sigma$, and
$\alpha$ be an element of $A$.
We call the subfunction
$$
 q_\alpha(x_1,\dots,x_{n-1}) \triangleq q(x_1,\dots,x_{n-1},\alpha).
$$
a {\it layer} of $q$.
It follows directly from \hyperref[p:11]{Proposition~11} and \hyperref[c:1]{Corollary~1} that

\proclaim{\wordCorollary~2}\label{c:2}
Let
$U_n$
be the set of partial $n$-qua\-si\-groups
$g:\Sigma^{n-1}\times\{a,b\}  \rightarrow \Sigma$
such that their layers
$g_\alpha$,
$\alpha\in\{a,b\}$
belong to a set
$W_{n-1}$ closed under equivalence.
Then
$|U_n|\leq (3w_{n-1}^2)/2^{n+1}$.
\endproclaim

\proclaim{\wordProposition~12 \rm (a representation of a decomposable $n$-qua\-si\-group
by the superposition of subfunctions)}\label{p:12}
Let $h$ and $g$ be an
$(n-m+1)$-
and $m$-qua\-si\-groups and
\begin{eqnarray}\nonumber
&f(x,\bar y, \bar z)\triangleq h(g(x,\bar y),\bar z),&
\\
h_0(x,\bar z)\triangleq f(x,\bar 0,\bar z),
&
g_0(x,\bar y)\triangleq f(x,\bar y,\bar 0),
&
\delta(x)\triangleq f(x,\bar 0,\bar 0)
                                \label{e:7}
\end{eqnarray}
where
$x\in\Sigma$,
$\bar y\in\Sigma^{m-1}$,
$\bar z\in\Sigma^{n-m}$.
Then
\begin{equation}
f(x,\bar y, \bar z)\equiv h_0(\delta^{-1}(g_0(x,\bar y)),\bar z).
                                \label{e:8}
\end{equation}
\endproclaim
\demo{Proof}
It follows from \hyperref[e:7]{(7)} that
$$
h_0(\cdot,\bar z)\equiv h(g(\cdot,\bar 0),\bar z),
\quad
g_0(x,\bar y)\equiv h(g(x,\bar y),\bar 0),
\quad
\delta^{-1}(\cdot)\equiv g^{\langle 1 \rangle}(h^{\langle 1 \rangle}(\cdot,\bar 0),\bar 0).
$$
Substituting these representations of
$h_0$,
$g_0$,
and
$\delta^{-1}$
to \hyperref[e:8]{(8)}, we can readily verify its validity.
\qed\enddemo

\proclaim{\wordProposition~13 \rm(on the number of the decomposable $n$-qua\-si\-groups)}\label{p:13}
For the number
$r^{\star}_n$
of the reduced decomposable $n$-qua\-si\-groups, it is true that
$$
r^{\star}_n\leq \sum_{m=2}^{n-1}{{n}\choose{m}}
v^{\star}_{n-m+1}v^{\star}_m.
$$
\endproclaim

\demo{Proof}
From \hyperref[p:12]{Proposition~12} we see that
a reduced decomposable $n$-qua\-si\-group can be represented
(maybe ambiguously) as a superposition
of reduced
$(n-m+1)$-
and $m$-qua\-si\-groups with
$m\in \{2,\dots,n-1\}$.
For every such $m$
the number of ways to split the set of arguments into two groups equals
${n}\choose{m}$;
and the numbers of ways to choose
$(n-m+1)$-
and $m$-qua\-si\-groups  equal respectively
$v^{\star}_{n-m+1}$
and
$v^{\star}_m$.
The order of arguments in each of the groups is not essential,
because
a reduced $m$-qua\-si\-group
goes into a reduced $m$-qua\-si\-group under a coordinate permutation.
\qed\enddemo

\section{Semilinear $n$-qua\-si\-groups}\label{s:4}

\demo{\wordDefinition}
An $n$-qua\-si\-group $f$ is called
{\it semilinear}
iff there are
$a,b\in \Sigma$
such that
$S_{a,b}(f)$
is a linear double-code.
An $n$-qua\-si\-group $f$ is called
{\it linear}
iff for all
$a,b\in \Sigma$,
$a\neq b$
the double-code
$S_{a,b} (f)$
is linear.
Fig.~3 gives illustrations of linear {\rm (a)},
semilinear {\rm (b)}, and nonsemilinear {\rm (c)} $3$-qua\-si\-groups.
\enddemo

\proclaim{\wordProposition~14}\label{p:14}
The reduced linear $n$-qua\-si\-group is unique.
\endproclaim

\demo{Proof}
The statement follows from \hyperref[p:4]{Proposition~4(e)}
and the fact that every $n$-qua\-si\-group $f$ is uniquely defined
by the double-MDS-codes
$S_{0,1}(f)$
and
$S_{0,2}(f)$.
\qed\enddemo

Denote by
$K_n\subseteq V_n$
the set of all semilinear $n$-qua\-si\-groups and by
$K_n(a,b)$
the set of semilinear $n$-qua\-si\-groups $f$
such that the double-code
$S_{a,b} (f)$
is linear. The validity of the following proposition is easily seen.

\proclaim{\wordProposition~15}\label{p:15}
For every different $a$, $b$, $c$ in
$\Sigma$
the intersection
$K_n(a,b) \cap K_n(a,c)$
is the set of all linear $n$-qua\-si\-groups.
\endproclaim

Using \hyperref[p:5]{Proposition~5(a)}, the following fact is easily proved by induction on $m$.

\proclaim{\wordProposition~16}\label{p:16}
Let $f$ be a semilinear $n$-qua\-si\-group. Then for every
$(a_1,\ldots,a_{m})\in \Sigma^{m}$
the function
$$
g(x_1,\ldots,x_{n-m}) \triangleq f(x_1,\ldots,x_{n-m},a_1,\ldots,a_{m})
$$
is a semilinear $(n-m)$-qua\-si\-group.
\endproclaim

\proclaim{\wordProposition~17}\label{p:17}
{\rm (a)}
 Equivalent $n$-qua\-si\-groups are
 or are not semilinear simultaneously.

{\rm (b)}
 If $f$ is a semilinear $n$-qua\-si\-group, then its inversions
 $f^{\langle i\rangle}$,
 $i\in [n]$
 also are semilinear $n$-qua\-si\-groups.
\endproclaim

\demo{Proof}
P.\,(a) follows from the fact that the set of linear double-codes is
closed under equivalence (\hyperref[p:4]{Proposition~4(a)}).

Let us prove p.\,(b). It can be checked directly that the semilinearity of $f$
is equivalent to the existence of
$a_0=a$,
$b_0=b$,
$a_1,\ldots,a_n$,
$b_1,\ldots,b_n$
such that
$a_i\neq b_i$
and
\begin{equation}
\bigoplus_{i=0}^n \chi_{\{a_i,b_i\}}(x_i) = 0
                                \label{e:9}
\end{equation}
for all
$x_0,x_1,\dots,x_n$
satisfying
$x_0=f(x_1,x_2,\dots,x_n)$.
Since the expression \hyperref[e:9]{(9)} is symmetric with respect to the choice
of the dependent variable, the claim is proved.
\qed\enddemo

\demo{\wordRemark}
The reduced linear $n$-qua\-si\-group $f$ can be represented in the form
$f(x_1,\dots,x_n) = x_1*\dots *x_n$
where
$(\Sigma,*)$
is a group isomorphic to
$Z_2\times Z_2$,
with the addition table

{\offinterlineskip
\def\extra{\omit\hfil&height2pt&\cr}
\def\tabrule{\noalign{\hrule}}
\def\strut{\vrule height10pt depth 4pt width0pt\relax}
$$
\vbox{
\halign {\strut\ \hfil$#$\hfil\ &\vrule#\tabskip=0pt plus 1fill\ &
\hfil \ $#$\ \hfil&\hfil \ $#$\ \hfil&
\hfil \ $#$\ \hfil& \hfil \ $#$\ \hfil
\tabskip=0pt\cr
*&&0&1&2&3\cr
\extra
\tabrule
\extra
0&&0&1&2&3\cr
1&&1&0&3&2\cr
2&&2&3&0&1\cr
3&&3&2&1&0\cr
}}
$$
}
\ifx\forTranslationNeeds\undefined
\else
\enddemo 
\fi

The following two lemmas were proved in \cite{Kro:2002ACCT:double-codes,Kro:2codes}.
The first concerns a representation of a nonprime double-MDS-code
by prime double-codes of smaller dimensions.
The second lemma, an essential corollary of the former,
connect the decomposability property of $n$-qua\-si\-group $q$ with
the nonprimality property
of $S_{c,d}(q)$.

\proclaim{\wordLemma~2 \rm(on decomposition of a double-MDS-code) \cite{Kro:2002ACCT:double-codes,Kro:2codes}}\label{l:2}
Let $S$ be a double-{\rm MDS}-code. Then there exists
$k=k(S)\in [n]$
such that

{\rm (a)} the characteristic function
 $\chi_S$
 can be represented as
\begin{equation}
 \chi_S(\bar x) \equiv \bigoplus_{j=1}^k \chi_{S_j}(\tilde x_j)
                                \label{e:10}
\end{equation}
 where
 $\tilde x_j=(x_{i_{j,1}},\dots ,x_{i_{j,n_j}})$
 are disjoint collections of variables from
 $\bar x$,
 $S_j\subset \Sigma^{n_j}$
 are prime double-{\rm MDS}-codes for $j\in [k];$
 the representation is unique up to substitution
 of double-{\rm MDS}-codes
 $S_j\backslash \Sigma^{n_j}$ for
 some double-{\rm MDS}-codes
 $S_j;$

{\rm (b)}
 $S$ is a union of
 $2^{k-1}$
 pairwise disjoint prime double-codes of equal cardinality$;$
 $\Sigma^n \backslash S$
 is a union of
 $2^{k-1}$
 pairwise disjoint prime double-codes of equal cardinality.
\endproclaim

\proclaim{\wordLemma~3 \rm(on the decomposability of $n$-qua\-si\-groups) \cite{Kro:2002ACCT:double-codes,Kro:2codes}}\label{l:3}
Let
$S\subset \Sigma^n$ be
a double-{\rm MDS}-code that satisfies {\rm \hyperref[e:10]{(10)}},
$c\neq d\in \Sigma$,
and  $q$ be
an $n$-qua\-si\-group such that
$S_{c,d}(q)=S$.
Then
\begin{equation}
 q(\bar x)  \equiv  q_0(q_1(\tilde x_1),\dots ,q_k(\tilde x_k))
                                \label{e:11}
\end{equation}
where
$q_j$, $j\in [k]$ are
$n_j$-qua\-si\-groups,
$q_0$ is
a semilinear $k$-qua\-si\-group, and the collections of variables
$\tilde x_j=(x_{i_{j,1}},\dots ,x_{i_{j,n_j}})$,
$j\in [k]$
and the numbers $k$,
$n_j$
are defined by \hyperref[l:2]{Lemma~{\rm 2}}.
\endproclaim

\proclaim{\wordCorollary~3}\label{c:3}
Let
$\{a,b,c,d\}=\Sigma$,
let $q$ be an $n$-qua\-si\-group, and let a partial $n$-qua\-si\-group
$g \triangleq q|_{\Sigma^{n-1}\times\{a,b\}}$
have more than two extensions. Then
$q\in R_n\cup K_n$.
\endproclaim

\demo{Proof}
It follows from definitions that
$C_a(f^{\langle n\rangle})=C(f_a)$
for an arbitrary $n$-qua\-si\-group $f$ and its inversion in $n$th argument
$f^{\langle n\rangle}$.
Let
$$
S \triangleq \Sigma^n\backslash (C(g_a)\cup C(g_b)).
$$
Then for every extension $f$ of the partial $n$-qua\-si\-group $g$ we see that
$$
 S = \Sigma^n\backslash (C(f_a)\cup C(f_b))
   = C(f_c)\cup C(f_d) = S_{c,d}(f^{\langle n\rangle}).
$$
By the hypothesis, the partial $n$-qua\-si\-group $g$ has more than two extensions $f$.
Each of the extensions is uniquely defined by its layer
$f_c$.
Hence the double-MDS-code $S$ includes more than two different MDS codes
$C(f_c)$.
By \hyperref[p:1]{Proposition~1},
the double-MDS-code
$S = S_{c,d}(q^{\langle n\rangle})$
consists of more than one prime double-code.
According to \hyperref[l:2]{Lemmas~2}  \hyperref[l:3]{and~3},
the number $k$ in \hyperref[e:11]{(11)} is not less than $2$.
If
$k<n$,
then \hyperref[e:11]{(11)} implies the decomposability of
$q^{\langle n\rangle}$;
if
$k=n$,
then \hyperref[e:10]{(10)} implies the semilinearity. So,
$q^{\langle n\rangle}\in K_n \cup R_n$;
then by \hyperref[p:9]{Propositions~9(b)} \hyperref[p:17]{and~17(b)} we get
$q\in K_n \cup R_n$.
\qed\enddemo

\section{On the number of $n$-qua\-si\-groups}\label{s:5}

In this section, we evaluate the number of the $n$-qua\-si\-groups of order $4$,
by establishing that the subclass of semilinear $n$-qua\-si\-groups is
asymptotically dominant.
We first calculate the number of the semilinear $n$-qua\-si\-groups.

\proclaim{\wordTheorem~1 \rm(on the number of the semilinear $n$-qua\-si\-groups)}\label{t:1}
$k^{\star}_n=3\cdot 2^{2^n-n-1} -2$ and $k_n=3^{n+1}\cdot 2^{2^n+1}-2^36^n$.
\endproclaim

\demo{Proof}
An arbitrary $n$-qua\-si\-group $f$ in
$K_n^{\star}(0,1)$
can be defined by firstly choosing the linear double-code
$S_{0,1}(f)$
and secondly, the MDS codes
$C_{0}(f)\subset S_{0,1}(f)$
and
$C_{2}(f)\subset \Sigma^n\backslash S_{0,1}(f)$.
A linear double-code can be chosen in
$2\cdot 3^n$
ways
(\hyperref[p:4]{Proposition~4(f)});
an MDS code,
in $2^{2^{n-1}}$
ways
(\hyperref[p:1]{Proposition~1}).
So,
$$
|K_n(0,1)|=2\cdot 3^n\cdot 2^{2^{n-1}}\cdot 2^{2^{n-1}} = 3^n\cdot 2^{2^n+1}.
$$
By \hyperref[c:1]{Corollary~1} we get
$|K^{\star}_n(0,1)|=2^{2^n-n-1}$
and, similarly,
$$
|K^{\star}_n(0,2)|=|K^{\star}_n(0,3)|=2^{2^n-n-1}.
$$
It follows from \hyperref[p:14]{Propositions~14} \hyperref[p:15]{and~15} that the pairwise intersections of
$K^{\star}_n(0,1)$,
$K^{\star}_n(0,2)$,
$K^{\star}_n(0,3)$
contain only one element.
Then, by the formula of inclusion and exclusion,
$$
k^{\star}_n = |K^{\star}_n(0,1)\cup K^{\star}_n(0,2)\cup
 K^{\star}_n(0,3)| = 3\cdot 2^{2^n-n-1} - 3 +1.
$$
By \hyperref[c:1]{Corollary~1}, we have
$k_n=4\cdot 6^nk^{\star}_n$.
\qed\enddemo

\demo{\wordRemark }
The lower bound
$v_n\geq 3^{n+1}\cdot 2^{2^n+1}-2^36^n$
was established in \cite{Kro:quas_l_b}.
\enddemo

As a result of a numerical experiment, we have the values:
\begin{equation}
v^{\star}_1 =1,
\quad
v^{\star}_2 =4,
\quad
v^{\star}_3 =64\quad \mbox{\cite{MulWeb}},
\quad
v^{\star}_4 =7132,
\quad
v^{\star}_5 =201538000.
                                \label{e:12}
\end{equation}

The following lemma shows that the existence of a semilinear layer in
a $n$-qua\-si\-group yields an arrangement of its structure.

\proclaim{\wordLemma~4 \rm(on a semilinear layer)}\label{l:4}
Let $q$ be an $n$-qua\-si\-group and there exists
$\alpha \in \Sigma$
such that
$q_\alpha\in K_{n-1}$.
Then
$q\in K_n\cup R_n$.
\endproclaim

\demo{Proof}
Assume that
$q_{\alpha}\in K_{n-1}$
for some
$\alpha\in \Sigma$
and thus the double-MDS-code
$S_{a,b}(q_{\alpha})$
is linear for some
$a,b\in \Sigma$.
Consider
$S_{a,b}(q)$;
we have
$S_{a,b}(q_{\alpha})= \TheLAYER{(S_{a,b}(q))}{n}{\alpha}$.
Then, by \hyperref[l:1]{Lemma~1}, there is
$\beta\in \Sigma$,
$\beta\neq \alpha$
such that
$$
S_{a,b}(q_{\beta}) = \TheLAYER{(S_{a,b}(q))}{n}{\beta}
  = \Sigma^{n-1} \backslash S_{a,b}(q_{\alpha}),
$$
i.\,e.,
the $(n-1)$-qua\-si\-group
$q_{\beta}$
is semilinear.

$(*)$
 We claim that the partial $n$-qua\-si\-group
$g \triangleq q|_{\Sigma^{n-1}\times\{\alpha,\beta\}}$
has two semilinear extensions.
Let
$\{a,b,c,d\}= \{\alpha, \gamma,\beta, \delta\} = \Sigma$
and
$\sigma \triangleq (ab)(cd)$
be a permutation of symbols of
$\Sigma$.
Then the function $f$ defined by the equalities
$$
f(x_1,\dots,x_{n-1},\alpha)  \triangleq   q(x_1,\dots,x_{n-1},\alpha),
 \
f(x_1,\dots,x_{n-1},\beta )  \triangleq   q(x_1,\dots,x_{n-1},\beta ),
$$
$$
f(x_1,\dots,x_{n-1},\gamma)  \triangleq  \sigma q(x_1,\dots,x_{n-1},\alpha),
 \
f(x_1,\dots,x_{n-1},\delta)  \triangleq  \sigma q(x_1,\dots,x_{n-1},\beta)
$$
is an extension of the partial $n$-qua\-si\-group $g$.
It is clear that
$S_{a,b}(f_\gamma) = S_{a,b}(f_\alpha)=S_{a,b}(q_\alpha)$;
therefore the double-codes
$\TheLAYER{(S_{a,b}(f))}{n}{\alpha} = \TheLAYER{(S_{a,b}(f))}{n}{\gamma}$
are linear; hence, by \hyperref[p:5]{Proposition~5(b)}, the double-code
$S_{a,b}(f)$
also is linear. So, the $n$-qua\-si\-groups $f$ and
$f'(\bar x) \triangleq f(x_1,\dots,x_{n-1},\tau(x_{n}))$
with
$\tau \triangleq (\gamma, \delta)$
satisfy $(*)$.

We note finally that either $q$ coincides with
one of $f$, $f'$, and thus
$q\in K_n$;
or $g$ has more than two extensions
($q$, $f$, $f'$), and
$q\in K_n\cup R_n$
by \hyperref[c:3]{Corollary~3}.
\qed\enddemo

\proclaim{\wordTheorem~2 \rm(on the number of the $n$-qua\-si\-groups)}\label{t:2}
If
$n\geq 5$,
then
$$
3^{n+1}2^{2^n +1}\leq v_n \leq (3^{n+1}+1)2^{2^n +1}.
$$
\endproclaim

\demo{Proof}
Let
$q\in V_n$;
consider the partial $n$-qua\-si\-group
$g_{\alpha, \beta} = q|_{\Sigma^{n-1}\times\{\alpha, \beta\} }$
for arbitrary
$\alpha, \beta \in \Sigma$.
If 
$g_{\alpha, \beta}$
has more than two extensions, then we have
$q\in K_n\cup R_n$, by \hyperref[c:3]{Corollary~3}.
If
$q_\alpha \in K_{n-1}$
or
$q_\beta \in K_{n-1}$,
then
$q \in K_n\cup R_n$, by \hyperref[l:4]{Lemma~4}.
Hence if
$q \not \in K_n\cup R_n$,
then for all
$\alpha,\beta \in \Sigma $
we have
$q_{\alpha},q_{\beta} \not \in K_{n-1}$
and the partial $n$-qua\-si\-group
$g_{\alpha, \beta}$
has two extensions.

Introduce the notation
$T_n \triangleq V_n\backslash K_n$
and
$W_n \triangleq T_n\backslash R_n$.
It follows from \hyperref[p:9]{Propositions~9(a)} \hyperref[p:17]{and~17(a)} that the sets
$T_n$
and
$W_n$
are closed under equivalence.
Then
$q\in W_n$
implies
$q_\alpha \in T_n$
for all
$\alpha \in \Sigma$
and, by \hyperref[c:2]{Corollary~2},
\begin{equation}
w_n \leq \frac {3t^2_{n-1}}{2^n}.
                                \label{e:13}
\end{equation}

$(*)$
We claim that the following three inequalities hold;
we prove them by induction on $n$.

\smallskip
{\rm (a)}
 $k^{\star}_n\leq v^{\star}_n\leq 2k^{\star}_n $
 whenever
 $n\geq 1$;

\smallskip
{\rm (b)}
 $t_{n}\leq 2^{2^n +1}$
 whenever
 $n\geq 5$;

\smallskip
{\rm (c)}
 $v_n \leq (3^{n+1}+1)2^{2^n +1} $
 whenever
 $n\geq 5$.

\smallskip
When
$n\leq 5$,
the conditions (a)--(c) are verified based on the exact values of
$k^{\star}_n$,
$v^{\star}_n$,
$v_n$,
$t_n=v_n-k_n$
(\hyperref[e:12]{(12)}, \hyperref[t:1]{Theorem~1}).
By the inductive assumption (a) holds for
$n\in [m]$,
and  (b), (c) hold for
$n=m\geq 5$.
Let us show the validity of (a)--(c) for
$n=m+1$.
From (a) and \hyperref[t:1]{Theorem~1}
with
$m\geq 5$,
$m-1>i>2$ we get
the following:
$$
{v^{\star}_{m-i+1}v^{\star}_i}\leq
4k^{\star}_{m-i+1}k^{\star}_i< 4\cdot9\cdot 2^{2^{m-i+1}+2^i-m-3}
<
4\cdot 3\cdot 2^{2^{m-1}-m-1} = v^\star_2k^{\star}_{m-1} \leq
{v^{\star}_{m-1}v^{\star}_2}.
$$
Since
$v^{\star}_2=4$,
from the estimate for
$r^\star_n$
(\hyperref[p:13]{Proposition~13}) we derive
$$
 r^\star_{m+1}
\leq \sum_{i=2}^{m}{{m+1}\choose{i}}
v^{\star}_{(m+1)-i+1}v^{\star}_i \leq
\sum_{i=2}^{m}{{m+1}\choose{i}} v^{\star}_{m}v^{\star}_2 <
2^{m+1}\cdot v^{\star}_m\cdot 4.
$$
Substituting (c) with
$n=m$,
we have
\begin{equation}
 r_{m+1} < 2^{m+3}(3^{m+1}+1)2^{2^m+1} <2^{2^{m+1}}.
                                \label{e:14}
\end{equation}
Moreover, from \hyperref[e:13]{(13)} and (b) with
$n=m$
we get the inequality
\begin{equation}
w_{m+1} \leq  \frac {3t^2_{m}}{2^{m+1}} \leq
\frac {3\cdot 2^{2^{m+1}+2}}{2^{m+1}} < 2^{2^{m+1}}.
                                \label{e:15}
\end{equation}
By the definitions of the sets
$T_m$
and
$W_m$
we have
$t_{m+1}\leq w_{m+1} + r_{m+1}$
and
$v_{m+1} = t_{m+1} + k_{m+1}$.
Then from the  inequalities \hyperref[e:14]{(14)} and \hyperref[e:15]{(15)} we derive 
(b) with
$n=m+1$,
and from \hyperref[t:1]{Theorem~1} and the inequality (b) we derive 
(a) and (c) with
$n=m+1$.
The claim $(*)$ is proved.

It remains to show the lower estimate for
$v_n$.
First we prove that
the following holds for $n\geq 4$:
\begin{equation}
 t^{\star}_n\geq t^{\star}_3v^{\star}_{n-2}.
                                \label{e:16}
\end{equation}
Let
$g\in T^{\star}_3$
and
$h\in V^{\star}_{n-2}$.
Then \hyperref[p:16]{Proposition~16} implies that the $n$-qua\-si\-group
$$
f(x_1,\dots,x_{n}) \triangleq h(g(x_{1},x_2,x_{3}),x_{4},\dots, x_{n})
$$
is not semilinear.
It is easy to check that different pairs of reduced
$(n-2)$-qua\-si\-group
$h$ and $3$-qua\-si\-group $g$ correspond to
different reduced $n$-qua\-si\-groups $f$.
The inequality \hyperref[e:16]{(16)} is proved.

From \hyperref[e:12]{(12)} and \hyperref[t:1]{Theorem~1} it follows that
$t^{\star}_3 = 18$.
Thus the inequality \hyperref[e:16]{(16)} and \hyperref[t:1]{Theorem~1} imply
$v^{\star}_n = k^{\star}_n + t^{\star}_n \geq 3^n2^{2^n-n-1}$
for
$n\geq 4$.
Then from \hyperref[c:1]{Corollary~1} we get the inequality
$v_n\geq 3^{n+1}2^{2^n +1}$
for
$n\geq 4$.
\qed\enddemo

It directly follows from \hyperref[t:2]{Theorem~2} and \hyperref[p:8]{Proposition~8} that
\proclaim{\wordCorollary~4 \rm(the asymptotic forms of the number of $n$-qua\-si\-groups
and the number of MDS codes)}\label{c:4}
Let
$m_n$
be the number of {\rm MDS}-codes in
$\Sigma^n$
and
$v_n$
be the number of $n$-qua\-si\-groups of order $4$. Then
$$
v_n = 3^{n+1}2^{2^n +1}(1+o(1)),\quad
m_n = 3^{n}2^{2^{n-1} +1}(1+o(1)).
$$
\endproclaim

\ifx\href\undefined \newcommand\href[2]{#2} \fi\ifx\url\undefined
  \newcommand\url[1]{\href{#1}{#1}} \fi\ifx\bbljan\undefined
  \newcommand\bbljan{Jan} \fi\ifx\bblfeb\undefined \newcommand\bblfeb{Feb}
  \fi\ifx\bblmar\undefined \newcommand\bblmar{March} \fi\ifx\bblapr\undefined
  \newcommand\bblapr{Apr} \fi\ifx\bblmay\undefined \newcommand\bblmay{May}
  \fi\ifx\bbljun\undefined \newcommand\bbljun{June} \fi\ifx\bbljul\undefined
  \newcommand\bbljul{July} \fi\ifx\bblaug\undefined \newcommand\bblaug{Aug}
  \fi\ifx\bblsep\undefined \newcommand\bblsep{Sep} \fi\ifx\bbloct\undefined
  \newcommand\bbloct{Oct} \fi\ifx\bblnov\undefined \newcommand\bblnov{Nov}
  \fi\ifx\bbldec\undefined \newcommand\bbldec{Dec} \fi

\end{document}